\documentclass[12pt]{amsart}
\newtheorem{theorem}{Theorem}[section]
\newtheorem{lemma}[theorem]{Lemma}
\newtheorem{proposition}[theorem]{Proposition}

\theoremstyle{definition}
\newtheorem{definition}[theorem]{Definition}
\newtheorem{example}[theorem]{Example}

\theoremstyle{remark}
\newtheorem{remark}[theorem]{Remark}
\numberwithin{equation}{section}
\newfont{\kh}{msbm10}
\newcommand{\C}{\mbox{\kh C}}
\begin{document}
\title {$n$-Homomorphisms}
\author{S. Hejazian}
\address{Department of Mathematics, Ferdowsi University, P. O. Box 1159, Mashhad 91775, Iran}
\email{hejazian@math.um.ac.ir}
\author{M. Mirzavaziri}
\address{Department of Mathematics, Ferdowsi University, P. O. Box 1159, Mashhad 91775, Iran}
\email{mirzavaziri@math.um.ac.ir}
\author{M. S. Moslehian}
\address{Department of Mathematics, Ferdowsi University, P. O. Box 1159, Mashhad 91775, Iran}
\email{moslehian@ferdowsi.um.ac.ir} 
\begin{abstract}
Let $\mathcal A$ and $\mathcal B$ be two (complex) algebras. A
linear map $\varphi:{\mathcal A}\to{\mathcal B}$ is called
$n$-homomorphism if $\varphi(a_{1}\ldots
a_{n})=\varphi(a_{1})\ldots\varphi(a_{n})$ for each
$a_{1},\ldots,a_{n}\in{\mathcal A}.$ In this paper, we
investigate $n$-homomorphisms and their relation to
homomorphisms. We characterize $n$-homomorphisms in terms of
homomorphisms under certain conditions. Some results related to
continuity and commutativity are given as well.
\end{abstract}
\maketitle
\section{Introduction}

Let $\mathcal A$ and $\mathcal B$ be two algebras. A linear
mapping $\varphi:{\mathcal A}\to{\mathcal B}$ is called an
$n$-homomorphism if $\varphi(a_{1}\ldots a_{n})=
\varphi(a_{1})\ldots\varphi(a_{n})$ for each
$a_{1},\ldots,a_{n}\in{\mathcal A}$. A $2$-homomorphism is then a
homomorphism, in the usual sense, between algebras.

For a homomorphism $\varphi:{\mathcal A}\to {\mathcal B}$ we can
see that $\varphi(a_{1}\ldots
a_{n})=\varphi(a_{1})\ldots\varphi(a_{n})$ for each
$a_{1},\ldots,a_{n}\in {\mathcal A}$ and for each $n$. The
problem is to consider the converse: Is every $n$-homomorphism a
homomorphism? The answer is negative in general. As a
counterexample, let $\omega$ be an $n$-root of unity and a
homomorphism $\psi: A \to B$. Then $\varphi : = \omega \psi$ is an
$n$-homomorphism which is not an $m$-homomorphism for any $2\leq
m\leq n-1$.

The paper consists of five sections. The second section is
devoted to examine the relationship between notions of
$n$-homomorphism and homomorphism. The third section concerns
$n$-homomorphisms on algebras generated by idempotents. In the
forth section, we investigate $n$-homomorphisms which preserve
commutativity under some conditions. A study of $n$-homomorphisms
on Banach algebras is given in the fifth section.

Throughout the paper, all Banach algebras are assumed to be over
the complex field ${\C}$.

\section{Relationship Between $n$-Homomorphisms and Homomorphisms}

We begin this section with a typical example as an exension of
that of introduced in the first section.

\begin{example} Let ${\mathcal A}$ be a unital algebra, $a_{0}$ be a
central element of ${\mathcal A}$ with $a_{0}^{n}=a_{0}$ for some
natural number $n$ and let $\theta:{\mathcal A}\to {\mathcal A}$
be a homomorphism. Define $\varphi:{\mathcal A}\to {\mathcal A}$
by $\varphi(a)=a_{0}\theta(a)$. Then we have
\begin{eqnarray*}
\varphi(a_{1}\ldots a_{n})&=&a_{0}\theta(a_{1}\ldots a_{n})\\
&=& a_{0}^{n}\theta(a_{1})\ldots\theta(a_{n})\\
&=& a_{0}\theta(a_{1})\ldots a_{0}\theta(a_{n})\\
&=& \varphi(a_{1})\ldots\varphi(a_{n}).
\end{eqnarray*}
Hence $\varphi$ is an $n$-homomorphism. In addition,
$a_{0}=\varphi(1_{\mathcal A})$ whenever $\theta$ is onto.
\end{example}

The above example gives us an $n$-homomorphism as a multiple of a
homomorphism. Indeed, if ${\mathcal A}$ has the identity
$1_{{\mathcal A}}$ then each $n$-homomorphism is of this form,
where $a_{0}=\varphi(1_{{\mathcal A}})$ as the following
proposition shows.

\begin{proposition} Let ${\mathcal A}$ be a unital algebra with identity $1_{{\mathcal A}}$, ${\mathcal B}$
be an algebra and $\varphi:{\mathcal A}\to {\mathcal B}$ be an
$n$-homomorphism. If $\psi:{\mathcal A}\to {\mathcal B}$ is
defined by $\psi(a)=(\varphi(1_{{\mathcal A}}))^{n-2}\varphi(a)$
then $\psi$ is a homomorphism and $\varphi(a)=\varphi(1_{\mathcal
A})\psi(a)$.
\end{proposition}

\begin{proof} We have
$$\varphi(1_{{\mathcal A}})=\varphi(1_{{\mathcal A}}^{n})=(\varphi(1_{{\mathcal A}}))^{n}.$$
and
\begin{eqnarray*}
\psi(ab)&=&(\varphi(1_{{\mathcal A}}))^{n-2}\varphi(ab)\\
&=&(\varphi(1_{{\mathcal A}}))^{n-2}\varphi(a 1_{{\mathcal A}}^{n-2}b)\\
&=&(\varphi(1_{{\mathcal
A}}))^{n-2}\varphi(a)(\varphi(1_{{\mathcal
A}}))^{n-2}\varphi(b)\\&=&\psi(a)\psi(b).
\end{eqnarray*}
It follows from $(\varphi(1_{{\mathcal A}}))^{n-1}\varphi(a)=
\varphi(1_{{\mathcal A}}^{n-1}a)=\varphi(a)$ that
$(\varphi(1_{{\mathcal A}}))^{n-1}$ is an identity for
$\varphi({\mathcal A})$. Thus
\begin{eqnarray*}
\varphi(1_{{\mathcal A}})\psi(a)&=&\varphi(1_{{\mathcal
A}})((\varphi(1_{{\mathcal A}}))^{n-2}\varphi(a))\\
&=&(\varphi(1_{{\mathcal A}}))^{n-1}\varphi(a)\\&=&\varphi(a).
\end{eqnarray*}
\end{proof}

Whence we characterized all $n$-homomorphisms on a unital algebra.
For a non-unital algebra ${\mathcal A}$ we use the unitization
and some other useful constructions. Recall that for an algebra
${\mathcal A}$, the linear space ${\mathcal A}_{1}={\mathcal
A}\oplus {\C}=\{(a,\alpha)\vert a\in {\mathcal A}, \alpha\in
{\C}\}$ equipped with the multiplication
$(a,\alpha)(b,\beta)=(ab+\alpha b+\beta a,\alpha\beta)$,
so-called the unitization of ${\mathcal A}$, is a unital algebra
with identity $(0,1)$ containing ${\mathcal A}$ as a two-sided
ideal.

Now we shall prove that each $n$-homomorphism is a multiple of a
homomorphism under some conditions.

\begin{definition} An algebra ${\mathcal A}$ is called a factorizable algebra
if for each $a\in {\mathcal A}$ there are $b,c\in {\mathcal A}$
such that $a=bc$.
\end{definition}

\begin{theorem} Let ${\mathcal A}$ and ${\mathcal B}$ be two factorizable algebras,
$lan({\mathcal B})=\{b\in {\mathcal B}; b{\mathcal B}=0\}=\{0\}$
and $\varphi:{\mathcal A}\to {\mathcal B}$ an onto
$n$-homomorphism. Then $\ker \varphi$ is a two-sided ideal of
${\mathcal A}$ and there is a unital algebra $\tilde{{\mathcal
B}}\supseteq {\mathcal B}$ and an $x\in \tilde{{\mathcal B}}$ with
$x^{n-1}=1_{\tilde{{\mathcal B}}}$ such that $\psi:{\mathcal A}\to
\tilde{{\mathcal B}}$ defined by $\psi(a)=x^{n-2}\varphi(a)$ is a
homomorphism.
\end{theorem}

\begin{proof} Suppose that $a\in \ker\varphi$ and $u\in {\mathcal A}$.
Since ${\mathcal A}$ is a factorizable algebra there are $u_1,
\ldots u_{n-1}\in {\mathcal A}$ such that $u=u_1\ldots u_{n-1}$.
Hence \[\varphi(au)=\varphi(au_1\ldots
u_{n-1})=\varphi(a)\varphi(u_1)\ldots \varphi(u_{n-1})=0.\]
Therefore $au\in \ker\varphi$. Similarly $ua\in \ker\varphi$.

Let $\tilde{{\mathcal
B}}=\{b_\circ+\beta_1x+\ldots+\beta_{n-2}x^{n-2}; b_\circ\in
{\mathcal B}_1, ~{\rm and }~\beta_1, \ldots, \beta_{n-2}\in {\bf
C}\}$ as a subset of the algebra ${\mathcal B}_1[x]$ of all
polynomials in $x$ with coefficients in the unitization
${\mathcal B}_1$ of ${\mathcal B}$. Using the ordinary
multiplication of polynomials, we define a multiplication on
$\tilde{{\mathcal B}}$ by $x^{n-1}=1$ and
$bx=\varphi(a_1)\varphi(a_2)$ where
$b=\varphi(a)=\varphi(a_1a_2)$ and $a=a_1a_2\in {\mathcal A}$. We
show that the multiplication is well-defined.

Let $b=d\in {\mathcal B}$ and $b=\varphi(a)=\varphi(a_1a_2),
d=\varphi(c)=\varphi(c_1c_2)$ with $a=a_1a_2, c=c_1c_2\in
{\mathcal A}$. Then we have $\varphi(a_1a_2)=\varphi(c_1c_2)$. So
$\varphi(a_1a_2)b_2\ldots b_n=\varphi(c_1c_2)b_2\ldots b_n$
for all $b_2\ldots b_n\in {\mathcal B}$.\\
Since $\varphi$ is onto there exist $u_2\ldots u_n\in {\mathcal
A}$ such that $\varphi(u_i)=b_i$. We can then write
\begin{eqnarray*}
&&\varphi(a_1)\varphi(a_2)\varphi(u_2)\ldots\varphi(u_{n-2})\varphi(u_{n-1}u_n)\\
&=&\varphi(a_1a_2u_2\ldots u_{n-1}u_n)\\
&=&\varphi(a_1a_2)\varphi(u_2)\ldots
\varphi(u_{n-1})\varphi(u_n)\\&=&
\varphi(c_1c_2)\varphi(u_2)\ldots
\varphi(u_{n-1})\varphi(u_n)\\
&=&\varphi(c_1c_2u_2\ldots
u_{n-1}u_n)\\&=&\varphi(c_1)\varphi(c_2)\varphi(u_2)
\ldots\varphi(u_{n-2})\varphi(u_{n-1}u_n).
\end{eqnarray*}
This implies that
$\varphi(a_1)\varphi(a_2)b=\varphi(c_1)\varphi(c_2)b$ for each
$b\in {\mathcal B}$, since ${\mathcal B}$ is a factorizable
algebra. Hence
$(\varphi(a_1)\varphi(a_2)-\varphi(c_1)\varphi(c_2)){\mathcal
B}=0$. Since $lan({\mathcal B})=\{0\}$ we conclude that
$\varphi(a_1)\varphi(a_2)=\varphi(c_1)\varphi(c_2)$. In
particular, $\varphi(a)\varphi(b)x^{n-2}=\varphi(ab)$ for all
$a,b\in {\mathcal A}$. Note that associativity of our
multiplication inherits from that of multiplication of
polynomials.

We can inductively prove that
$\varphi(a_1)\ldots\varphi(a_m)x^{n-m}= \varphi(a_1\ldots a_m)$
for all $m\geq 2$.

To show this, suppose that it holds for $m\geq 2$ and $a_{m+1}\in
{\mathcal A}$. Then
\begin{eqnarray*}
&&\varphi(a_1)\ldots\varphi(a_{m-1})\varphi(a_m)\varphi(a_{m+1})x^{n-m-1}\\
&=&\varphi(a_1)\ldots\varphi(a_{m-1})\varphi(a_m)\varphi(a_{m+1})x^{n-(m+1)}x^{n-1}\\
&=&\varphi(a_1)\ldots\varphi(a_{m-1})(\varphi(a_m)\varphi(a_{m+1})x^{n-2})x^{n-m}\\
&=&\varphi(a_1)\ldots\varphi(a_{m-1})\varphi(a_ma_{m+1})x^{n-m}\\&=&\varphi(a_1\ldots
a_{m-1}a_ma_{m+1}).
\end{eqnarray*}

Now define $\tilde{\varphi}:{\mathcal A}_{1}\to \tilde{{\mathcal
B}}$ by $\tilde{\varphi}(a,\alpha)=\varphi(a)+\alpha x$ for each
$(a,\alpha)\in {\mathcal A}_{1}$. Then for each
$(a_{1},\alpha_{1}), \ldots,(a_{n},\alpha_{n})\in {\mathcal
A}_{1}$ we have
\[\tilde{\varphi}(\prod_{i=1}^{n}(a_{i},\alpha_{i}))=
\tilde{\varphi}(\sum \alpha_{j_{1}}\ldots\alpha_{j_{k}}a_{i_{1}}
\ldots a_{i_{l}}),\] where the summation is taken over all
$i_{1},\ldots,i_{l},j_{1}, \ldots,j_{k}$ with $i_{1}<\ldots
<i_{l}, j_{1}<\ldots <j_{k}, 0\leq k,l\leq n,
\{i_{1},\ldots,i_{l}\}\cap\{j_{1},\ldots,j_{k}\}=\emptyset,
\{i_{1},\ldots,i_{l}\}
\cup\{j_{1},\ldots,j_{k}\}=\{1,\ldots,n\}.$ Thus if $\varphi()$
denotes $1\in {\C}$ then we can write

\begin{eqnarray*}
\tilde{\varphi}(\prod_{i=1}^{n}(a_{i},\alpha_{i}))
&=&\sum \alpha_{j_{1}}\ldots\alpha_{j_{k}}\varphi(a_{i_{1}}\ldots a_{i_{l}})\\
&=&\sum \alpha_{j_{1}}\ldots\alpha_{j_{k}}\varphi(a_{i_{1}})
\ldots\varphi(a_{i_{l}})x^{k}\\
&=&\prod_{i=1}^{n}(\varphi(a_{i})+\alpha_{i}x)=\prod_{i=1}^{n}
\tilde{\varphi}(a_{i},\alpha_{i}).
\end{eqnarray*}

This shows that $\tilde{\varphi}$ is an $n$-homomorphism on
${\mathcal A}_{1}$. Now Proposition 2.3 implies that
$\tilde{\psi}: {\mathcal A}_{1}\to\tilde{{\mathcal B}}$ defined by
$\tilde{\psi}(a,\alpha)=(\tilde{\varphi}(1_{{\mathcal
A}_{1}}))^{n-2}
\tilde{\varphi}(a,\alpha)=(\tilde{\varphi}(0,1))^{n-2}(\varphi(a)+\alpha
x) =x^{n-2}(\varphi(a)+\alpha x)$ is a homomorphism on ${\mathcal
A}_{1}$. Thus $\psi: {\mathcal A}\to \tilde{{\mathcal B}}$
defined by $\psi(a)=x^{n-2}\varphi(a)$ is a homomorphism on
${\mathcal A}$.
\end{proof}

\begin{proposition} Let ${\mathcal A}$ and ${\mathcal B}$ be two algebras,
$\psi:{\mathcal A}\to {\mathcal B}$ be a homomorphism and
$\tilde{{\mathcal B}}$ be a unital algebra with $\tilde{{\mathcal
B}}\supseteq {\mathcal B}$. Then for each central element $x\in
\tilde{{\mathcal B}}$ with $x^{n-1}=1_{\tilde{{\mathcal B}}}$,
the mapping $\varphi:{\mathcal A}\to \tilde{{\mathcal B}}$
defined by $\varphi(a)=x\psi(a)$ is an $n$-homomorphism and
$\ker\varphi$ is a two-sided ideal of ${\mathcal A}$.
\end{proposition}

\begin{proof} We have
\begin{eqnarray*}
\varphi(a_1\ldots a_n)&=&x\psi(a_1\ldots
a_n)\\&=&x^n\psi(a_1)\ldots
\psi(a_n)\\&=&x\psi(a_1)x\psi(a_2)\ldots
x\psi(a_n)\\&=&\varphi(a_1)\ldots \varphi(a_n).
\end{eqnarray*}
Moreover if $a\in \ker\varphi$ and $b\in {\mathcal A}$ then
$\varphi(ab)=x\psi(ab)=x\psi(a)\psi(b)=\varphi(a)\psi(b)=0$ and
similarly $\varphi(ba)=0$. Thus $\ker\varphi$ is a two-sided ideal
of ${\mathcal A}$.
\end{proof}

\begin{example} In general, the kernel of an $n$-homomorphism
may not be an ideal. As an example, take algebra ${\mathcal A}$
of all $3$ by $3$ matrices having $0$ on and below the diagonal.
In this algebra product of any $3$ elements is equal to $0$ so
any linear map from ${\mathcal A}$ into itself is a
$3$-homomorphism but its kernel does not need to be an ideal.
\end{example}

\section{Algebras Generated By Idempotents}

Let ${\mathcal A}$ and ${\mathcal B}$ be two algebras. A linear
map $\varphi:{\mathcal A}\to {\mathcal B}$ is called associative
if $\varphi(a)\varphi(bc)=\varphi(ab)\varphi(c)$ for all $a, b,
c\in {\mathcal A}$. For example every mapping preserving zero
products is associative (see Lemma 2.1 of \cite{CKLW}.

\begin{proposition} Let ${\mathcal A}$ be linearly generated by its
idempotents and $\varphi:{\mathcal A}\to {\mathcal B}$ be an
$n$-homomorphism. Then
$\varphi(a)\varphi(bc)=\varphi(ab)\varphi(c)$ for all $a, b, c\in
{\mathcal A}$.
\end{proposition}

\begin{proof} Let $a, c\in {\mathcal A}$ and $e\in {\mathcal A}$ be an
idempotent. Then
\begin{eqnarray*}
\varphi(a)\varphi(ec)&=&\varphi(a)\varphi(e^{n-1}c)\\
&=&\varphi(a)\varphi(e)\ldots\varphi(e)\varphi(c)\\&=&\varphi(ae^{n-1})\varphi(c)\\
&=&\varphi(ae)\varphi(c).
\end{eqnarray*}
Since $\varphi$ is linear and the linear span of idempotent
elements of ${\mathcal A}$ is ${\mathcal A}$ itself, we obtain
$\varphi(a)\varphi(bc)=\varphi(ab) \varphi(c)$ for all $a, b,
c\in {\mathcal A}$.
\end{proof}

The following lemma states some interesting properties of associative mappings.

\begin{lemma} If ${\mathcal A}$ is unital and $\varphi:{\mathcal A}\to {\mathcal B}$ is associative then\\
{\rm (i)} $\varphi(1)$ belongs to the center of image of $\varphi$;\\
{\rm (ii)} For all $a, b\in {\mathcal A}, \varphi(1)\varphi(ab)=\varphi(a)\varphi(b)$;\\
{\rm (iii)} $\varphi$ preserves commutativity;\\
{\rm (iv)} If ${\mathcal B}$ is unital, and $\varphi(1)$ is
invertible or the algebra $<\varphi({\mathcal A})>$ generated by
$\varphi({\mathcal A})$ has an identity, then
$\psi(a)=\varphi(1)^{-1}\varphi(a)$ is a homomorphism. {\rm ( see
Lemma 2.1 of \cite{CKLW})}
\end{lemma}

\begin{proof} We merely prove (iv), the rest is easily
followed from the definition.

If $\varphi(1)$ is invertible in $<\varphi({\mathcal A})>$ and
$\psi(a)=\varphi(1)^{-1}\varphi(a)$, then
\begin{eqnarray*}
\psi(ab)&=&\varphi(1)^{-1}\varphi(ab)\\&=&
\varphi(1)^{-1}\varphi(1)^{-1}\varphi(a)\varphi(b)\\&=&
\varphi(1)^{-1}\varphi(a)\varphi(1)^{-1}\varphi(b)\\&=&\psi(a)\psi(b).
\end{eqnarray*}
If the algebra $<\varphi({\mathcal A})>$ has an identity, say
$u$, then $u=\displaystyle{\sum_{i=1}^n}\varphi(a_{i_1})
\ldots\varphi(a_{i_n})$. Hence
\[u=uu=(\displaystyle{\sum_{i=1}^n}
\varphi(a_{i_1})\ldots\varphi(a_{i_n}))^2=\varphi(1)v=v\varphi(1).\]
Therefore $\varphi(1)$ is invertible in $<\varphi({\mathcal A})>$.
\end{proof}

\begin{remark} The statements of Lemma 3.2 hold for every
unital $C^*$-algebra of real rank zero or $W^*$-algebra
${\mathcal A}$ and every norm-continuous associative mapping
$\varphi$ from ${\mathcal A}$ into a normed algebra.
\end{remark}

The following result is known for mappings preserving zero
products; cf. Corollary 2.8 of \cite{CKLW}. We however reproved
it for the sake of ellegance.

\begin{proposition} Suppose that $H$ and $K$ are infinite
dimensional complex Hilbert spaces and $\varphi:B(H)\to B(K)$ is
a bijective $n$-homomorphism, then there exist $\alpha\in {\C}$
and an invertible bounded linear operator $S$ from $K$ onto $H$
such that $\varphi(T)=\alpha S^{-1}TS$ for all $T\in B(H)$.
\end{proposition}

\begin{proof} It is known that every bounded linear operator
on an infinite dimensional complex Hilbert space is a sum of at
most five idempotent \cite{P-T}. Hence $B(H)$ is generated by its
idempotents. Since $\varphi(1)\neq 0$ belongs to the center the
factor $B(K)$, it is a scalar $\alpha\in {\C}$. Then Lemma 3.2
follows that $\psi=\alpha^{-1}\varphi$ is an isomorphism from
$B(H)$ onto $B(K)$. Applying Theorem 4 of \cite{ARN}, there
exists an invertible bounded linear operator $S:K\to H$ such that
$\psi(T)=S^{-1}TS$ for all $T\in B(H)$.
\end{proof}

\section{Commutativity}

Recall that an algebra ${\mathcal A}$ is called semiprime if
$a{\mathcal A}a=\{0\}$ implies that $a=0$ for each $a\in {\mathcal
A}$. The known techniques similar to \cite{BRE} yield the
following:

\begin{lemma} Let ${\mathcal A}$ be a semiprime algebra and let $a\in {\mathcal
A}$. If $[[a,x],x]=0$ for all $x\in {\mathcal A}$, then $[a,x]=0$
for all $x\in {\mathcal A}$.
\end{lemma}

\begin{proof} For all $x,y\in {\mathcal A}$,
\[0=[[a,x+y],x+y]=[[a,x]+[a,y], x+y]=[[a,x],y]+[[a,y],x] ~~~(1)\]
Taking $yx$ for $y$ in (1) we have
\begin{eqnarray*}
0&=&[[a,x],yx]+[[a,yx],x]\\&=&[[a,x],y]x+[[a,y]x+
y[a,x],x]\\&=&[[a,x],y]x+[[a,y],x]x+[y,x][a,x].
\end{eqnarray*}
Thus for all $x,y\in {\mathcal A}$,
\[ [y,x][a,x]=0~~~(2).\]
Replacing y by az in (2) we obtain
\begin{eqnarray*}
0&=&[az,x][a,x]\\&=&([a,x]z+a[z,x])[a,x]\\&=&[a,x]z[a,x]+a[z,x][a,x]\\&=&[a,x]z[a,x].
\end{eqnarray*}
Hence $[a,x]{\mathcal A}[a,x]=\{0\}. ~~~(3)$

Since ${\mathcal A}$ is semiprime (3) implies that $[a,x]=0$.
\end{proof}

\begin{theorem} Suppose that ${\mathcal A}$ and ${\mathcal B}$ are two
algebras, ${\mathcal B}$ is semiprime and $\varphi:{\mathcal
A}\to {\mathcal B}$ is a surjective $3$-homomorphism. If
${\mathcal A}$ is commutative, then so is ${\mathcal B}$.
\end{theorem}

\begin{proof} For each $a, b\in {\mathcal A}$, $[[a,b],b]=[0,b]=0$. So
that
\[[[\varphi(a),\varphi(b)],\varphi(b)] = \varphi([[a,b],b]) = 0.\]
Applying lemma 4.1, $[\varphi(a),\varphi(b)]=0$. Hence ${\mathcal
B}$ is commutative.
\end{proof}

\section{$n$-Homomorphisms on Banach Algebras}

Recall that the second dual ${\mathcal A}^{**}$ of a Banach
algebra ${\mathcal A}$ equipped with the first Arens product is a
Banach algebra. The first Arens product is indeed characterized
as the unique extension to ${\mathcal A}^{**} \times {\mathcal
A}^{**}$ of the mapping $(a,b) \mapsto ab$
from ${\mathcal A} \times {\mathcal A}$ into ${\mathcal A}$ with the following properties :\\
(i) for each $G \in {\mathcal A}^{**}$, the mapping $F \mapsto FG$ is weak*-continuous on ${\mathcal A}^{**}$;\\
(ii) for each $a \in {\mathcal A}$, the mapping $G \mapsto aG $ is weak*-continuous on ${\mathcal A}^{**}$.\\
The second Arens product can be defined in a similar way. If the
first and the second Arens products coincide on ${\mathcal
A}^{**}$,
then ${\mathcal A}$ is called regular.\\
We identify ${\mathcal A}$ with its image under the canonical
embedding $i : {\mathcal A} \longrightarrow {\mathcal A}^{**}$.

\begin{theorem} Suppose that ${\mathcal A}$ and ${\mathcal B}$ are two Banach algebras
and $\varphi:{\mathcal A}\to {\mathcal B}$ is a continuous
$n$-homomorphism. Then the second adjoint $\varphi^{**}:{\mathcal
A}^{**}\to {\mathcal B}^{**}$ of $\varphi$
is also an $n$-homomorphism.\\
If, in addition, ${\mathcal A}$ is Arens regular and has a
bounded approximate identity, then $\phi$ is a certain multiple
of a homomorphism.
\end{theorem}

\begin{proof} Let $F_1, \ldots, F_n\in {\mathcal A}^{**}$.
By Goldstine' theorem (cf. \cite{D-S}), there are nets $(a^1_i),
\ldots, (a^n_j)$ in ${\mathcal A}$ such that \[{\rm weak}^*{\rm
-}\lim _i a^1_i=F_1, \ldots, {\rm weak}^*{\rm -}\lim _j
a^n_j=F_n.\]
Since $\varphi^{**}$ is ${\rm weak}^*{\rm
-}$continuous we have
\begin{eqnarray*}
\varphi^{**}(F_1 \ldots F_n)&=&\varphi^{**} ({\rm weak}^*{\rm
-}\lim _i\ldots {\rm weak}^*{\rm -}\lim _j a^1_i\ldots
a^n_j)\\&=& {\rm weak}^*{\rm -}\lim _i\ldots {\rm weak}^*{\rm
-}\lim _j \varphi^{**}(a^1_i\ldots a^n_j)\\&=& {\rm weak}^*{\rm
-}\lim _i\ldots {\rm weak}^*{\rm -}\lim _j \varphi^(a^1_i\ldots
a^n_j)\\&=&{\rm weak}^*{\rm -}\lim _i \ldots {\rm weak}^*{\rm
-}\lim _j (\varphi(a^1_i)\ldots \varphi(a^n_j))\\&=& {\rm
weak}^*{\rm -}\lim _i\varphi(a^1_i)\ldots {\rm weak}^*{\rm -}\lim
_j \varphi(a^n_j)\\&=& {\rm weak}^*{\rm -}\lim
_i\varphi^{**}(a^1_i) \ldots {\rm weak}^*{\rm -}\lim _j
\varphi^{**}(a^n_j)\\&=& \varphi^{**}(F_1)\ldots\varphi^{**}(F_n).
\end{eqnarray*}

If ${\mathcal A}$ is Arens regular and has a bounded approximate
identity , it follows from and Proposition 28.7 of \cite{B-D} that
${\mathcal A}^{**}$ has an identity. By proposition 2.3 there
exists a homomorphism $\psi:{\mathcal A}^{**}\to {\mathcal
B}^{**}$ such that $\varphi(a)=\varphi^{**}|_{\mathcal
A}(a)=\varphi^{**}(1_{{\mathcal A}^{**}})\psi(a)$ for all $a\in
{\mathcal A}$.
\end{proof}

\begin{remark} A computational proof similar to that of
Theorem 6.1 of \cite{C-Y} may be used in extending
$n$-homomorphisms to the second duals.
\end{remark}

Now suppose that $\varphi$ is a non-zero $3$-homomorphism from a
unital algebra ${\mathcal A}$ to ${\C}$. Then $\varphi(1)= 1$ or
$-1$. Hence either $\varphi$ or $-\varphi$ is a character on
${\mathcal A}$. If ${\mathcal A}$ is a Banach algebra, then
$\varphi$ is automatically continuous; cf. Theorem 16.3 of
\cite{B-D}. It may however happen that a $3$-homomorphism is not
continuous:

\begin{example} Let ${\mathcal A}$ be the algebra of all $3$ by $3$ matrices
having $0$ on and below the diagonal and ${\mathcal B}$ be the
algebra of all ${\mathcal A}$-valued continuous functions from
$[0,1]$ into ${\mathcal A}$ with sup norm. Then ${\mathcal B}$ is
an infinite dimensional Banach algebra ${\mathcal B}$ in which the
product of any three elements is $0$. Since $\mathcal B$ is
infinite dimensional there are linear discontinuous maps (as
discontinuous $3$-homomorphisms) from ${\mathcal B}$ into itself.
\end{example}

\begin{theorem} Let ${\mathcal A}$ be a $W^*$-algebra and ${\mathcal B}$ a $C^*$-algebra.
If $\varphi:{\mathcal A}\to {\mathcal B}$ is a weakly-norm
continuous 3-homomorphism preserving the involution, then $\Vert
\varphi\Vert\leq 1$.
\end{theorem}

\begin{proof} The closed unit ball of ${\mathcal A}$ is compact in weak
topology. By the Krein-Milman theorem this convex set is the
closed convex hull of its extreme points. On the other hand, the
extreme points of the closed unit ball of ${\mathcal A}$ are the
partial isometries $x$ such that $(1-xx^*){\mathcal
A}(1-x^*x)=\{0\}$, cf. Problem 107 of \cite{HAL}, and Theorem
I.10.2 of \cite{TAK}. Since
$\varphi(xx^*x)=\varphi(x)\varphi(x)^*\varphi(x)$, the mapping
$\varphi$ preserves the partial isometries. Since every partial
isometry $x$ has norm $\Vert x\Vert\leq 1$, we conclude that
$$\Vert \varphi(\displaystyle{\sum_{i=1}^n}\lambda_ix_i)\Vert=
\Vert\displaystyle{\sum_{i=1}^n}\lambda_i\varphi(x_i)\Vert\leq
\displaystyle{\sum_{i=1}^n}\lambda_i\Vert\varphi(x_i)\Vert\leq 1$$
where $x_1, \ldots, x_n$ are partial isometries,
$\lambda_1\ldots \lambda_n>0$ and $\displaystyle{\sum_{i=1}^n}\lambda_i=1$.\\
It follows from weak continuity of $\varphi$ that
$\Vert\varphi\Vert\leq 1$.
\end{proof}

{\bf Question.} Is every $*$-preserving $n$-homomorphism between
$C^*$-algebras continuous?

{\bf Aknowledgment.} The authors would like to thank Professor
Krzysztof Jarosz for providing the Examples 2.6 and 5.3.

\bibliographystyle{amsplain}

\end{document}